# A construction of magic 24-cells

Park, Donghwi

Seoul National University

Abstract

I found a novel class of magic square analogue, magic 24-cell. The problem is to assign the consecutive numbers 1 through 24 to the vertices in a graph, which is composed of 24 octahedra and 24 vertices, to make the sum of the numbers of each octahedron the same.

It is known that there are facet-magic and face-magic labelings of tesseract. However, because of 24-cell contains triangle, face-magic labeling to assign different labels to each vertex is impossible.

So I tried to make a cell-magic labeling of 24-cell. Linear combination of three binary labeling and one ternary labeling gives 64 different magic labelings of 24-cell.

Due to similarity in the number of vertices between 5x5 magic square and magic 24-cell, It might be possible to calculate the number of magic 24-cell.

Further analysis would be need to determine the number of magic 24-cell.

## 0. Introduction

There are many magic square varieties which assign face-magic labeling to the vertices of planar graphs. Notable examples are Jisugwimundo[1], Yeonhwando[2] and etc.

Therefore, it is of interest to construct face-magic labelling of regular polyhedra However, because of regular deltahedra are consist by triangles, face-magic labelings of regular deltahedra are impossible. Furthermore, magic normal labeling of regular dodecahedron is impossible due to parity. However, there are face-magic labelings of the cube.[Fig. 1]

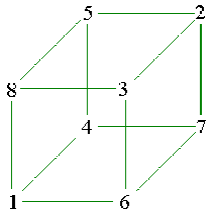

Fig 1. Face-magic cube.

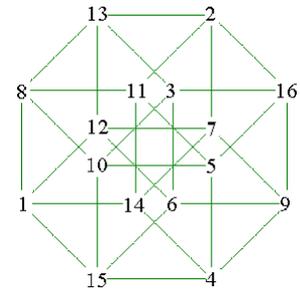

Fig 2. Tesseract-Cube-Square [3]

By this way, we can generalize face-magic labelings of regular polyhedra to facet-magic or face-magic labelings of regular polychora.

It is known that there are facet-magic and face-magic labelings of tesseract.[3] [Fig. 2]
However, face-magic labelings of regular 4-polytopes which are constructed by odd-gons are impossible due to parity. Facet-magic labelings of regular 4-polytopes which are constructed by tetrahedra are also impossible.

The remaining cases are facet-magic labelings of 24-cell or 120-cell.

The 24-cell is the convex regular 4-polytope which is composed of 24 octahedral cells. The 120-cell is the convex regular 4-polytope which is composed of 120 dodecahedral cells. Due to difficulty on constructing facet-magic labelings of 120-cell, I focused on to construct facet-magic labelings of 24-cell.

1. Coordinate representation and indexing of vertices

The 24 vertices of 24-cell can be described as the 24 coordinate permutations of: $(\pm 1, \pm 1, 0, 0)$. Then, I assigned indices 1 to 24 to the vertices of 24-cell in descending order.
Then, I assigned indices on 24 octahedra and assign six vertices from each octahedron.
8 octahedra was obtained by permuting the integer coordinates: $(\pm 1, 0, 0, 0)$.
16 octahedra was obtained by element-wise multiplication of $(\pm 1, \pm 1, \pm 1, \pm 1)$ and 6 coordinate permutations of $(1,1,0,0)$.

2. Parities of numbers in magic 24-cell.

In magic 24-cell, the magic sum of octahedra should be 75. It means every vertex sum of octahedra should be odd. So each octahedron should contain odd number of odd numbers.

If we assign 0 or 1 to 24-cell, each octahedron should contain one, three or five 1s and 0s. There are 2704156(=c(24,12) ) binary labelings on 24-cell which contains 12 ones and 12 zeros.

Brute-force search proves there are 256 binary labelings which contains odd numbers of odd number in each polyhedra.

In 192 binary labelings, some octahedra contains five or one odd numbers in each octahedron.
In 64 binary labelings, each octahedra contains three 0s and three 1s. These 64 labeling will be more useful to construct magic 24-cell.

In these 256 binary labelings, parities of opposite vertices were different.

Then, I made a ternary labelings on 24-cell. 24-cell can be decomposed into three 16-cells and color its vertices. [4]. So we can assign 0, 1 and 2 to each vertices and each triangle contains 0, 1 and 2 one times.

## 3. Linear combination.

We can pick three magic-binary labelling and one magic-ternary labelling. If their superimposition yields all possible 24 combinations of entries, all possible linear combinations give magic 24-cell.

If we set three binary (0 or 1) labelings as $B1$, $B2$ and $B3$. and ternary (0,1 or 2) labelings as $T$ and their superimposition yields all possible 24 combinations of entries. Then $12*B1+6*B2+3*B3+T$, $12*B1+6*B2+2*T+B3$, $12*B1+4*T+2*B2+B3$ and $8T+4*B1+2*B2+B3$ give different magic 24-cell.

When I fixed ternary labeling and its weights, there are 3072 superimpositions which yield all possible 24 combinations of entries.

However we can change the weight of ternary labelings in 4 ways and we can change the labels of ternary labelings by 6 ways, There are 73728 ways of magic 24-cell made by this linear combination which include rotation and reflection.
However, full symmetry group of the 24-cell has order 1152. So If we disregard rotation and reflection, we found 64 magic 24-cell.

## 4. Conclusion.

Linear combination of a few-color labelings gives an approach to construct the magic 24-cell.

This result does not exclude any other magic 24-cell labelings. Due to similarity in the number of vertices between 5x5 magic square and magic 24-cell, it might be possible to calculate the exact numbers of magic 24-cell.

Further analysis would be need to determine the total numbers of magic 24-cell.

Linear combinations of a few-color labelings might give an approach to construct the facet-magic 120-cell. 5-color labelings, ternary labelings, binary or quaternary labelings would be useful to assign 600 different values on the 600 vertices of 120-cell.

My implementation is available on Github.
https://github.com/gwahak/mathematics

## 5. Note

Jisugwimundo, also known as Hexagonal tortoise problem (HTP), is a magic square variety which was invented by medieval Korean Mathematician and minister Suk-Jung Choi (1646-1715).[1] In Jisugwimundo, numbers are arranged into vertices on a graph that is composed of hexagons. Also every hexagonal sum must be same.

Yeonhwando is a magic square variety which was first appeared at Yang Hui Suanfa. It contains octagons and squares.[2] In Yeonhwando, every octagonal sum must be same and every square sum must be same. Also every octagonal sum must be the twice of the square sum.